\documentclass[12pt,a4paper]{article}
\usepackage[latin2]{inputenc}
\usepackage{amsmath}
\usepackage{amsfonts}
\usepackage{amssymb}
\usepackage{graphicx}
\usepackage[left=2.00cm, right=2.00cm, top=2.50cm, bottom=2.50cm]{geometry}
\usepackage[T1]{fontenc}

\def\Q{{\mathbb Q}}
\def\Z{{\mathbb Z}}

\newtheorem{lemma}{Lemma}
\newtheorem{theorem}[lemma]{Theorem}

\title{
Power integral bases in a family of octic fields
}
\author{
Istv\'an Ga\'al\\
{\small University of Debrecen, Mathematical Institute} \\
{\small H--4002 Debrecen Pf.400., Hungary,} \\
{\small e--mail: gaal.istvan@unideb.hu},
}

\begin{document}
\baselineskip=17pt

\maketitle
\thispagestyle{empty}

\renewcommand{\thefootnote}{\arabic{footnote}}
\setcounter{footnote}{0}

\vspace{0.5cm}


\begin{abstract}
Several recent results prove the monogenity of some polynomials.
In these cases the root of the polynomial generates a power integral basis in the
number field generated by the root. A straightforward question is whether
such a  number field admits other generators of power integral bases?
We have investigated this problem in some previous papers 
 and here we extend this research to
a family of octic polynomials, following a recent result of 
L. Jones \cite{jones2}.
\end{abstract}

\section{Introduction}

A number field $K$ of degree $n$ with ring of integers $\Z_K$ 
is called {\it monogenic} (cf. \cite{book}) if there 
exists $\xi\in \Z_K$ such that $(1,\xi,\ldots,\xi^{n-1})$ is an integral basis, 
called power integral basis. 
We call $\xi$ the generator of this power integral basis.
$\alpha,\beta\in\Z_K$ are called {\it equivalent},
if $\alpha+\beta\in\Z$ or $\alpha-\beta\in\Z_K$. 
Obviously, $\alpha$ generates a power integral basis in $K$
if and only if any $\beta$, equivalent to $\alpha$, does.
As it is known, any algebraic number field admits up to equivalence only
finitely many generators of power integral bases.

An irreducible polynomial $f(x)\in\Z[x]$ is called {\it monogenic}, 
if a root $\xi$ of $f(x)$ generates a power integral basis in $K=\Q(\xi)$.
If $f(x)$ is monogenic, then $K$ is monogenic, but the converse is
not true.

For $\alpha\in\Z_K$ (generating $K$ over $\Q$) we call the module index 
\[
I(\alpha)=(\Z_K:\Z[\alpha])
\]
the {\it index} of $\alpha$. $\alpha$ generates a power integral basis in $K$
if and only if $I(\alpha)=1$. If $\alpha^{(i)}\; (1\le i\le n)$ 
are the conjugates of 
$\alpha$ in $K$  of degree $n$ with discriminant $D_K$, then
\[
I(\alpha)=\frac{1}{\sqrt{|D_K|}}\prod_{1\le i<j\le n}|\alpha^{(i)}-\alpha^{(j)}|.
\]
For more details concerning monogenity 
and power integral bases cf. \cite{book}.

In some recent papers we investigated number fields 
generated by a root of a  monogenic polynomial and made
calculations to figure out, whether these fields admit any additional
generators of power integral bases. We refer to 
\cite{g21} for some sextic trinomials,
\cite{g23b} for pure sextic fields,
\cite{g24t} for pure octic fields,
\cite{g24c} for certain quartic trinomials
and 
\cite{g24b} for some quartic polynomials with given Galois groups.

L. Jones \cite{jones1} (see also \cite{jones2}) gave conditions
for the monogenity of certain even octic polynomials
of type $x^8+ax^6+bx^4+ax^2+1$.
In the present paper we extend our calculations to this type of polynomials.
Among others we prove the existence of a non-trivial generator 
of power integral basis.

The octic field, generated by a root of the above polynomial 
is a quadratic extension of a quartic field.
It is an interesting point of our arguments, that this octic field
 can also be considered as a quartic extension of a quadratic field,
which makes it much easier to deal with. 

All tools used in our 
calculations are optimized to this special case, in order to 
make our calculations more efficient.

\section{The octic polynomial}

Let
\begin{equation}
f(x)=x^8+ax^6+bx^4+ax^2+1
\label{f}
\end{equation}
with $a,b\in \Z$. Set
$W_1 = b + 2 - 2a, W_2 = b + 2 + 2a, W_3 = a^2 - 4b + 8$.
L. Jones \cite{jones2} proved:
\begin{theorem}
If $W_1W_2W_3$ is square free and 
\begin{equation}
(a \mod 4,\; b\mod 4)\in\{(1,3),(3,1),(3,3)\}
\label{ab}
\end{equation}
then the polynomial $f(x)$ in (\ref{f}) is monogenic.
\label{j2}
\end{theorem}

Assume $f(x)$ is monogenic, not necessarily satisfying (\ref{ab}).
We wonder how many generators of power integral bases the number field $K$
has, generated by a root $\alpha$ of $f(x)$.

Set 
\begin{equation}
g(x)=x^4+ax^3+bx^2+ax+1.
\label{g}
\end{equation}
Obviously, a root $\beta$ of $g(x)$ satisfies $\alpha^2=\beta$.
Therefore the octic number field $K=\Q(\alpha)$ is a quadratic extension of
the quartic field $L=\Q(\beta)$:
\[
\Q\overset{4}{\subset} L=\Q(\beta) \overset{2}{\subset} K=\Q(\alpha).
\]
Unfortunately there exist feasible algorithms for solving 
index form equations, that is for determining generators of power integral bases,
only for a restricted class of number fields. Apart from low degree fields,
like cubic and quartic fields, there exist such algorithms only for 
some higher degree fields with special structure. These are e.g. sextic fields
with a quadratic subfield and octic fields with a quadratic subfield
(cf. \cite{book}). Above we have an octic field with a quartic
subfield, but using the reciprocal structure of $f(x)$ and $g(x)$
we can help this problem.

If
\[
\beta^4+a\beta^3+b\beta^2+a\beta+1=0,
\]
then
\[
\beta^2+a\beta+b+\frac{a}{\beta}+\frac{1}{\beta^2}=0,
\]
hence
\[
\left(\beta^2+\frac{1}{\beta^2}\right)+a\left(\beta+\frac{1}{\beta}\right)+b=0,
\]
\[
\left(\beta+\frac{1}{\beta}\right)^2+a\left(\beta+\frac{1}{\beta}\right)
+b-2=0.
\]
This yields, that 
\begin{equation}
\delta=\beta+\frac{1}{\beta}
\label{ddd}
\end{equation}
satisfies the quadratic equation
\begin{equation}
\delta^2+a\delta+(b-2)=0.
\label{delta}
\end{equation}
Consequently the number field $M=\Q(\delta)$ is a quadratic
subfield of $K$:
\[
\Q\overset{2}{\subset} M=\Q(\delta) \overset{4}{\subset} K=\Q(\alpha).
\]
By (\ref{ddd}) we have
\[
\beta^2-\delta\beta+1=0,
\]
therefore
\[
\alpha^4-\delta\alpha^2+1=0.
\]
This means that 
\begin{equation}
h(x)=x^4-\delta x^2+1
\label{h}
\end{equation}
is the relative defining polynomial
of $\alpha$ over $M$, and this is what we need for our procedure.

Note that in some previous papers (cf. \cite{gp}) we have developed
an algorithm for the complete resolution of index form equations
in octic fields with a quadratic subfield. This takes quite a long CPU time, since one has to solve a unit equation in the octic field. Also,
it can take long to calculate the fundamental units of that field, which
is a necessary input data for the calculations.

Therefore, if we would like to have an overall picture about the 
generators of power integral bases of our octic fields, 
we have to restrict ourselves to the calculation of the so called
"small solutions", that means we calculate all generators of 
power integral bases having coefficients, say $\le 10^{200}$
in absolute value in a given integral basis. 
Since the generators of power integral bases usually have very small
coefficients, such an algorithm determines all generators of
power integral bases with a very high probability. Moreover it certainly
indicates, if a number field, generated by a root of a monogenic polynomial,
has also other generators of power integral bases, in addition to 
the root of the polynomial.

Moreover, as we shall see in the following, a crucial point in this 
algorithm is the resolution of a relative quartic Thue equation
over the quadratic subfield. The fast algorithm \cite{relthue} 
for determining "small"
solutions of quartic relative Thue equations over quadratic fields
is only efficient, if the quadratic subfield is complex.
Therefore in our calculations we assume
\begin{equation}
a^2-4b+8<0,
\label{deltacomplex}
\end{equation} 
which guarantees by (\ref{delta}), that $M$ is a complex quadratic subfield.

On the other hand, we shall not restrict ourselves to those
monogenic polynomials $f(x)$, satisfying all conditions of Theorem \ref{j2}. 
We shall run the parameters $a,b$ in certain 
regions and consider all irreducible polynomials $f(x)$ that are monogenic.
The only condition we keep is that $W_3=a^2-4b+8$ is square-free,
in order to fix the basis element of $M$ and to make our arguments
simpler. Note, that we made calculations also for non-squarefree $W_3$,
and had comletely the same experiences, including also the non-trivial 
generator of power integral bases (cf. Theorm \ref{nnn}).

\section{Integral basis}

We have $M=\Q(\delta)$, with
\begin{equation}
\delta=\frac{-a+\sqrt{a^2-4b+8}}{2}.
\label{deltak}
\end{equation}
According to the above arguments, we assume 
$W_3=a^2-4b+8<0$ is square-free. To keep usual notation we
set $m=W_3$. This number can only be square-free if 
$a\equiv \; \pm 1\; (\bmod \; 4)$, whence $m\equiv \; 1\; (\bmod \; 4)$,
therefore the integral basis of the complex quadratic field
$M=\Q(\sqrt{m})$ is $(1,\omega)$,
where
\[
\omega=\frac{1+\sqrt{m}}{2}.
\]

We shall make calculations for monogenic polynomials $f(x)$.
In this case a root $\alpha$ of
$f(x)$ generates a power integral basis in $K=\Q(\alpha)$.

We shall use the following statements of \cite{grs}
which are certainly well known:
\begin{theorem}\mbox{}\\
A. If $K$ is monogenic, then $K$ is also relative monogenic over the subfield
$M$.\\
B. All generators of power integral bases of $K$ are of the form
\[
\gamma=X_0+\varepsilon \gamma_0,
\]
where $X_0\in\Z_M$, $\varepsilon$ is a unit in $M$ and $\gamma_0$
generates a relative power integral basis of $K$ over $M$.
\label{th2}
\end{theorem}

If 
\[
(1,\alpha,\alpha^2,\alpha^3,\alpha^4,\alpha^5,\alpha^6,\alpha^7)
\]
is an integral basis of $K$, then by the
 first part of the theorem 
\[
(1,\alpha,\alpha^2,\alpha^3)
\]
is a relative integral basis of $K$ over $M$, that is
any $\gamma\in\Z_K$ can be written in the form
\begin{equation}
\gamma=C+X\alpha+ Y\alpha^2+Z\alpha^3,
\label{XYZ}
\end{equation}
where
\[
C=c_1+\omega c_2, X=x_1+\omega x_2,Y=y_1+\omega y_2,Z=z_1+\omega z_2\in\Z_M
\]
with $c_1,c_2,x_1,x_2,y_1,y_2,z_1,z_2\in\Z$.

\section{A quartic relative Thue equations}

A consequence of \cite{gp} is the following

\begin{lemma}
Let $x^4+a_1x^3+a_2x^2+a_3x+a_4 \in \Z_M[x]$ be the relative defining polynomial
of $\alpha$ over $M$. Let
\begin{eqnarray*}
F(u,v)&=&u^3-a_2u^2v+(a_1a_3-4a_4)uv^2+(4a_2a_4-a_3^2-a_1^2a_4)v^3, \\
Q_1(x,y,z)&=&x^2 -xya_1 +y^2a_2+xz(a_1^2-2a_2)+yz(a_3-a_1a_2) \nonumber \\
                &&   +z^2(-a_1a_3+a_2^2+a_4) = u, \\
Q_2(x,y,z)&=&y^2-xz-a_1yz+z^2a_2 = v.
\end{eqnarray*}
If $\gamma_0$, represented in the form (\ref{XYZ})
 generates a relative power integral basis of $K$ over $M$
(that is the relative index of $I_{K/M}(\gamma_0)=(Z_K:\Z_M[\gamma_0])$
is equal to 1), 
then there exist $U,V\in\Z_M$ such that
\begin{eqnarray}
N_{M/\Q}(F(U,V))&=&\pm 1, \label{fuv}\\
Q_1(X,Y,Z)&=&U, \label{q1}\\
Q_2(X,Y,Z)&=&V. \label{q2}
\end{eqnarray}
\label{rellemma}
\end{lemma}

In our case by (\ref{h}) we have $a_1=0,a_2=-\delta,a_3=0,a_4=1$, hence
\[
F(u,v)=(u-2v)(u+2v)(u+\delta v).
\]
If $\gamma_0$ generates a relative power integral basis of $K$ over $M$,
then in view of the above Lemma, together with the $X,Y,Z\in\Z_M$
appearing in its representation (\ref{XYZ}) there exist $U,V\in \Z_M$
with
\[
N_{M/\Q}(F(U,V))=\pm 1.
\]
If $F(U,V)$ is a unit in $M$, then $U-2V,U+2V,U+\delta V$
are also units in $M$. Therefore 
\[
U-2V=\varepsilon_1,U+2V=\varepsilon_2
\]
and
\[
4V=\varepsilon_2-\varepsilon_1.
\]
In a complex quadratic field each unit is of absolute value 1,
and $V\in\Z_M$ is of absolute value 0 or $\ge 1$. Since the right side
is of absolute value $\le 2$, the above equation implies $V=0$.
As a consequence, $U$ is a unit in $M$. Then we have
\[
Q_1(X,Y,Z)=\varepsilon,\; Q_2(X,Y,Z)=0
\]
with a unit $\varepsilon\in M$. Following the arguments of \cite{gp}
we construct
\[
Q_0(X,Y,Z)=UQ_2(X,Y,Z)-VQ_1(X,Y,Z)=0,
\]
whence
\[
Q_2(X,Y,Z)=Y^2-XZ-\delta Z^2=0.
\]
$X_0=1,Y_0=0,Z_0=0$ is a non-trivial solution of $Q_2(X,Y,Z)=0$.
Using an argument of L. J. Mordell \cite{mordell} we parametrize
$X,Y,Z$ with $R,P,Q\in M$:
\begin{eqnarray}
X&=&RX_0\nonumber\\
Y&=&RY_0+P\label{pq}\\
Z&=&RZ_0+Q\nonumber
\end{eqnarray}
Substituting this representation of $X,Y,Z$ into $Q_2(X,Y,Z)=0$
we obtain 
\[
RQ=P^2-\delta Q^2.
\]
We multiply by $Q$ the equations in (\ref{pq}) and replace $RQ$
by $P^2-\delta Q^2$, then
\begin{equation}
\begin{array}{ccccc}
kX&=&P^2&& -\delta\; Q^2,\\
kY&=&&PQ, &  \\
kZ&=&&&Q^2,
\end{array}
\label{kpqpq}
\end{equation}
with a $k\in M$.
Further, applying the arguments of \cite{gp},
in (\ref{kpqpq}) we can replace the parameters $k,P,Q\in M$
by integer parameters in $\Z_M$, and it follows from the form of the
above coefficent matrix of $P^2,PQ,Q^2$  (and the property of
$\gamma_0$ being a generator of a relative power integral basis)
that $k$
is a unit in $M$. Finally, we substitute the representation
(\ref{kpqpq}) into $Q_1(X,Y,Z)=U$ and then we obtain
\begin{equation}
F(P,Q)=P^4-\delta P^2Q^2+Q^4=\varepsilon,
\label{FPQ}
\end{equation}
with a unit $\varepsilon$. 
This is a quartic relative Thue equation over the quadratic subfield $M$.
As $F(x,1)$ is just the relative defining polynomial of $\alpha$ 
over $M$, the equation can be written in the form
\begin{equation}
N_{K/M}(P-\alpha Q)=\varepsilon.
\label{NFPQ}
\end{equation}

\section{Solving the quartic relative Thue equation}

$M$ is a complex quadratic field, therefore the conjugate of
any $\nu\in M$ is its complex conjugate $\overline{\nu}$.
Denote by $\alpha^{(1)},\alpha^{(2)},\alpha^{(3)},\alpha^{(4)}$
the relative conjugates of $\alpha\in K$ over $M$,
corresponding to $\omega$
(these are the roots of $h(x)$ in (\ref{h})), then 
$\overline{\alpha^{(1)}},\overline{\alpha^{(2)}},
\overline{\alpha^{(3)}}, \overline{\alpha^{(4)}}$
are the relative conjugates of $\alpha$
over $M$ corresponding to $\overline{\omega}$.
Set $P=p_1+\omega p_2,Q=q_1+\omega q_2$ with $p_1,p_2,q_1,q_2\in\Z$.

Let $P,Q\in\Z_M$ be an arbitraty but fixed solution of (\ref{NFPQ}).
The unit $\varepsilon$ in (\ref{NFPQ}) is of absolute value 1, hence
using $\beta=P-\alpha Q$ (\ref{NFPQ}) implies
\begin{equation}
|\beta^{(1)}\beta^{(2)}\beta^{(3)}\beta^{(4)}|=1.
\label{bbb}
\end{equation}
Denote by $i_0$ the conjugate with
\[
|\beta^{(i_0)}|=\min_{1\le j\le 4} |\beta^{(j)}|.
\]
(We have to perform all calculations for all possible values of $i_0$.)
Then by (\ref{bbb}) we have 
$|\beta^{(i_0)}|\le 1$, whence
\begin{equation}
|P|\le   |\beta^{(i_0)}|+ \overline{|\alpha|} |Q|  \le 1+\overline{|\alpha|} |Q|,
\label{beta}
\end{equation}
where we denote by $\overline{|\alpha|}$ the size of $\alpha$,
that is the maximum absolute value of its conjugates.

Our purpose is to determine $c_2,x_1,x_2,y_1,y_2,z_1,z_2$ in (\ref{XYZ})
with absolute value $\le S=10^{200}$. This implies $|X|=|x_1+\omega x_2|
\le (1+|\omega|)S$ and similarly $|Z|\le (1+|\omega|)S$.
The representation (\ref{kpqpq})  of $Z$ implies
\[
|Q|\le \sqrt{|Z|}\le \sqrt{(1+|\omega|)S}.
\]
The representation of $X$ implies 
\[
|P|^2\le |X|+|\delta| |Q|^2\le |X|+|\delta| |Z|
\]
whence
\[
|P|\le \sqrt{(1+|\omega|)(1+|\delta|)S}.
\]
therefore 
\begin{equation}
\max (|P|,|Q|)\le  \max (|X|,|Y|) \le  \sqrt{(1+|\omega|)(1+|\delta|)S}.
\label{PQom}
\end{equation}
We have 
\begin{equation}
|p_1|= 
\frac{|\overline{\omega}P-\omega\overline{P}|}{|\overline{\omega}-\omega|}
\le \frac{2|\omega||P|}{|\overline{\omega}-\omega|},\;\;
|p_2|= \frac{|P-\overline{P}|}{|\omega-\overline{\omega}|} 
\le \frac{2|P|}{|\omega-\overline{\omega}|}
\label{p1p2}
\end{equation}
and similarly
\[
|q_1|\le  \frac{2|\omega||Q|}{|\overline{\omega}-\omega|},\;\;
|q_2|\le  \frac{2|Q|}{|\omega-\overline{\omega}|}.
\]
These imply 
\begin{equation}
A=\max(|p_1|,|p_2|,|q_1|,|q_2|)\le \frac{2|\omega|}{|\omega-\overline{\omega}|}  
\max(|P|,|Q|)\le   
\frac{2|\omega|}{|\omega-\overline{\omega}|} \sqrt{(1+|\omega|)(1+|\delta|)S},
\label{apq}
\end{equation}
that is we have to determine the solutions of (\ref{NFPQ})
until this bound. Note that for $S=10^{200}$ this bound is of magnitude 
$10^{100}$.

Further, together with (\ref{beta}) we have
\begin{equation}
A\le \frac{2|\omega|}{|\omega-\overline{\omega}|}\max(|P|,|Q|)
\le c_1 |Q|,
\label{Akorl}
\end{equation}
with
\[
c_1=\frac{2|\omega|}{|\omega-\overline{\omega}|}(1+\overline{|\alpha|}).
\]
If $|Q|\ge 10$ then for $1\le j\le 4,j\ne i_0$ this yields 
\begin{equation}
|\beta^{(j)}|\ge |\beta^{(j)}-\beta^{(i_0)}|-|\beta^{(i_0)}|
\ge |\alpha^{(j)}-\alpha^{(i_0)}| |Q|-1\ge (|\alpha^{(j)}-\alpha^{(i_0)}|-0.1) |Q|.
\label{bk}
\end{equation}
In our calculations we have to check all possible $q_1,q_2$ with $|Q|< 10$
separately.
(\ref{bbb}) and (\ref{bk}) imply
\begin{equation}
|\beta^{(i_0)}|=\frac{1}{\underset{\underset{j\ne i_0}{1\le j\le 4}}{\prod}|\beta^{(j)}|}
\ge \frac{1}{\underset{\underset{j\ne i_0}{1\le j\le 4}}{\prod}(|\alpha^{(j)}-\alpha^{(i_0)}|-0.1)} |Q|^{-3}\le c_{2,i_0} A^{-3},
\label{bbee}
\end{equation}
with
\[
c_{2,i_0}=\frac{c_1^3}{\underset{\underset{j\ne i_0}{1\le j\le 4}}{\prod}(|\alpha^{(j)}-\alpha^{(i_0)}|-0.1)}
\]
(depending on $i_0$).

\section{Reduction}

We apply a reduction procedure to reduce the bound in (\ref{Akorl}),
using inequality (\ref{bbee}), that is
\begin{equation}
|p_1+\omega p_2 -\alpha^{(i_0)} q_1 -\omega \alpha^{(i_0)} q_2|\le c_2 A^{-3}.
\label{redineq}
\end{equation}

We follow the arguments of \cite{relthue}.
Let $H$ be a large constant to be determined appropriately
(for a practical choice of $H$ see later).
Consider the lattice generated by the columns of the matrix
\[
\left(
\begin{array}{cccc}
1&0&0&0\\
0&1&0&0\\
0&0&1&0\\
0&0&0&1\\
   H&H\Re(\omega)&H\Re(-\alpha^{(i_0)})&H\Re(-\alpha^{(i_0)}\omega)\\
   0&H\Im(\omega)&H\Im(-\alpha^{(i_0)})&H\Im(-\alpha^{(i_0)}\omega)\\
\end{array}
\right).
\]

\vspace{0.5cm}

\begin{lemma} (cf. \cite{relthue}, or Lemma 5.3 of \cite{book})
Denote by $\ell_1$ the first vector of the LLL reduced basis of this lattice.
If $A\le A_0$ and $H$ is large enough to have
\begin{equation}
|\ell_1|\ge \sqrt{40}\cdot A_0,
\label{e1}
\end{equation}
then
\begin{equation}
A\le \left(\frac{c_{2,i_0}\cdot H}{A_0}\right)^{1/3}.
\label{redA}
\end{equation}
\label{redlemma}
\end{lemma}
Note that this procedure must be performed for all possible values of $i_0$.

We start with the upper bound $A_0$ in (\ref{apq}). For a certain $A_0$
usually $A_0^2$, $10\cdot A_0^2$ or $100 \cdot A_0^2$ 
is a suitable choice for $H$.
We have to make $H$ so large that (\ref{e1}) is satisfied.  
In view of (\ref{redA}) the new bound for $A$ will be of magnitude $A_0^{1/3}$
in the first reduction steps. The the following steps the reduction is not so 
fast anymore, but in about 8-10 steps the original bound of magnitude
$10^{100}$ is reduced to about 10.  A tipical sequence is the following:

\[
\begin{array}{|c|c|c|c|}
\hline
step & A_0 &H & new \; A_0 \\ \hline
1 & 10^{100} & 10^{202} & 9.1198\cdot 10^{33}\\ \hline
2 & 9.1198\cdot 10^{33} &  8.3172\cdot 10^{69} & 8.8440\cdot 10^{11}\\ \hline
3 &  8.8440\cdot 10^{11} & 7.8217\cdot 10^{25} & 87540.0136\\ \hline 
4 & 87540.0136 & 7.6632\cdot 10{11} & 187.9568 \\ \hline
5 & 187.9568 & 3.5327\cdot 10^6 & 24.2479 \\ \hline
6 & 24.2479 & 58796.3577 & 12.2522 \\ \hline
7 & 12.2522 & 15011.6532 & 9.7587 \\ \hline
\end{array}
\]
The reduction process is very fast, it usually only takes a few seconds.
For a constant $H$ of magnitude $10^{200}$ we have to use multiply precision
arithmetics with about 250 digits.

\section{Determining $P$ and $Q$}

We return to the quartic relative Thue equation (\ref{FPQ}), that is
\begin{equation}
P^4-\delta P^2Q^2+Q^4-\varepsilon=0.
\label{F0}
\end{equation}
Let $A_R$ be the reduced bound for $A$ obtained in the previous section.

By $m\equiv 1\; (\bmod \; 4)$, $|Q|<10$ yields 
$|q_1+\frac{1+\sqrt{m}}{2} q_2|<10$ whence $|q_2|<20/\sqrt{|m|}$ and 
$|q_1|<10+|q_2|/2<10+10/\sqrt{|m|}$.

Set $S_1=10+10/\sqrt{|m|}$, 
$S_2=20/\sqrt{|m|}$. Let $A_1=\max(A_R,S_1),A_2=\max(A_R,S_2)$.

We let $q_1$ run up to $|q_1|\leq A_1$ and $q_2$ run up to $|q_2|\leq A_2$.
For each pair $(q_1,q_2)$ we calculate $Q=q_1+\omega q_2$, substitute 
it into (\ref{F0}), and for all possible unit $\varepsilon\in M$ we solve
the quartic polynomial equation (\ref{F0}) for the complex number $P$.
Having the real and complex parts of $P$ we can determine $p_1,p_2$
with $P=p_1+\omega p_2$ (similarly as in (\ref{p1p2}))
and check if these values of $p_1,p_2$ are integers.

Having $P$ and $Q$ we can determine $X,Y,Z$ from (\ref{kpqpq}). Recall that
$k$ in  (\ref{kpqpq}) is a unit of $M$. Therefore all generators of
relative power integral bases of  $K$ over $M$ are of the form
$C+\varepsilon (\alpha X+\alpha^2 Y + \alpha^3Z)$ with arbitrary
$C\in\Z_M$ and arbitrary unit $\varepsilon\in M$.

\section{Determining generators of power integral bases of $K$}

For all possible  $X,Y,Z$ as calculated above,
we set $\gamma_0=\alpha X+\alpha^2 Y + \alpha^3 Z$.
In view of Theorem \ref{th2} all generators of power integral bases
are of the form 
\begin{equation}
\gamma=c_1+\omega c_2+\varepsilon \gamma_0,
\label{ggg}
\end{equation}
with $c_1,c_2\in\Z_M$, $\varepsilon$ is a unit in $M$.
In order to determine all non-equivalent generators of power integral
bases of $M$ we have to determine $\varepsilon$ and $c_2$ so that
$I(\gamma)=1$. For this purpose we shall use the following consequence
of Proposition 1 of \cite{grs}. 
Here we denote by $\gamma^{(1,j)}$ the conjugates of
$\gamma$ corresponding to $\alpha^{(j)}$ and by $\gamma^{(2,j)}$ the conjugates of $\gamma$ corresponding to $\overline{\alpha^{(j)}}$
for $1\le j\le 4$.

\begin{lemma}
\[
I(\gamma)=I_{K/M}(\gamma)\cdot J(\gamma)
\]
where
\[
I_{K/M}(\gamma)=\frac{1}{\sqrt{|N_{M/\Q}(D_{K/M})|}}
\prod_{i=1}^2\prod_{1\le j_1\le j_2\le 4}|\gamma^{(i,j_1)}-\gamma^{(i,j_2)}|
\]
is the relative index of $\alpha$ and
\[
J(\gamma)=\frac{1}{|D_M|^2}
\prod_{j_1=1}^4\prod_{j_2=1}^4 |\gamma^{(1,j_1)}-\gamma^{(2,j_2)}|.
\]
\label{relindex}
\end{lemma}

In view of Lemma \ref{rellemma} we calculated $\gamma_0$ to have
relative index 1. Any $\gamma$ of type (\ref{ggg}) is relative equivalent
to $\gamma_0$, that is their relative indices are equal.
Therefore we have to determine $\varepsilon$ and $c_2$ using $J(\gamma)=1$.
For all of the few possible units $\varepsilon$ of $M$, 
we calculate $J(\gamma)$. The equation 
\[
\prod_{j_1=1}^4\prod_{j_2=1}^4 
(\gamma^{(1,j_1)}-\gamma^{(2,j_2)})\pm D_M^2=0
\]
is a polynomial equation
with rational integer coefficients of degree 16.
In order to determine the possible values (if any) of $c_2\in\Z$,
corresponding to $\varepsilon$, we have to determine the integer roots
in $c_2$ of this polynomial. Note that $D_M=a^2-4b+8$.

\section{Results of our calculations}

Our routines were written in Maple.
We made calculations for several pairs $(a,b)$, such that
the polynomial $f(x)$ is irreducible, monogenic
 and $m$ is square-free. 
These pairs seldom satisfied the conditions of Theorm \ref{j2}.
Note that we also made calculations in cases when $m$ is not
squarefree (then $m=m_0\cdot m_1^2$ with square-free $m_1$
and $\omega$ can be either $(1+\sqrt{m_1})/2$ or $\sqrt{m_1}$)
and we had similar experiences.

The table below summarizes generators of power integral bases of 
$K$, represented in the form 
\[
\gamma=(c_1+\omega c_2)+(x_1+\omega x_2)\alpha+
(y_1+\omega y_2)\alpha^2+(z_1+\omega z_2)\alpha^3.
\]
We let $(a,b)$ run in $-25\le a\le 25, 2\le b\le 25$
and took those pairs $(a,b)$ 
for which $f(x)$ is irreducible, monogenic
and $m$ is square-free. In these 51 examples it took 526 
seconds (using an average PC) to calculate all generators of 
power integral bases with coefficients $\le 10^{200}$
in absolute value. We list $(a,b,m)$ and then the coefficients
$[c_2,x_1,x_2,y_1,y_2,z_1,z_2]$ of generators of power integral 
bases. We omit the trivial [0, 1, 0, 0, 0, 0, 0].

\[
\begin{array}{lllll}
(-9,23,-3),&
[0, 4, 1, 0, 0, -1, 0]&&&
\\
(-7,15,-3),&
[0, 3, 1, 0, 0, -1, 0]&&&
\\
(-7,19,-19),&
[0, 3, 1, 0, 0, -1, 0]&&&
\\
(-7,23,-35),&
[0, 3, 1, 0, 0, -1, 0]&&&
\\
(-5,9,-3),&
 [0, 2, 1, 0, 0, -1, 0],&
 [1, -2, 1, 1, -1, 1, -1], &
[1, 2, -1, 1, -1, -1, 1]&
\\
(-5,10,-7),&
 [0, 2, 1, 0, 0, -1, 0]&&&
\\
(-5,11,-11),&
[0, 2, 1, 0, 0, -1, 0]&&&
\\
(-5,14,-23),&
 [0, 2, 1, 0, 0, -1, 0]&&&
\\
(-5,18,-39),&
 [0, 2, 1, 0, 0, -1, 0]&&&
\\
(-5,19,-43),&
 [0, 2, 1, 0, 0, -1, 0]&&&
\\
(-5,21,-51),&
 [0, 2, 1, 0, 0, -1, 0]&&&
\\
(-5,22,-55),&
 [0, 2, 1, 0, 0, -1, 0]&&&
\\
(-5,23,-59),&
 [0, 2, 1, 0, 0, -1, 0]&&&
\\
(-5,25,-67),&
 [0, 2, 1, 0, 0, -1, 0]&&&
\\
(-3,7,-11),& 
 [0, 1, 1, 0, 0, -1, 0]&&&
\\
(-3,15,-43),&
 [0, 1, 1, 0, 0, -1, 0]
\\
(-1,3,-3),&
 [0, 0, 1, 0, 0, -1, 0],& 
[0, 1, -1, 0, 0, 0, 1],&
[0, 1, 0, 0, 0, -1, 1],&
[0, 1, -1, 0, 0, 0, 0],   \\&&
 [0, 0, 1, 0, 0, 0, 0]&&
\\
(-1,6,-15),&
 [0, 0, 1, 0, 0, -1, 0]&&&
\\
(-1,7,-19),&
 [0, 0, 1, 0, 0, -1, 0]&&&
\\
(-1,10,-31),&
 [0, 0, 1, 0, 0, -1, 0]&&&
\\
(-1,11,-35),&
 [0, 0, 1, 0, 0, -1, 0]&&&
\\
(-1,13,-43),&
 [0, 0, 1, 0, 0, -1, 0]&&&
\\
(-1,15,-51),&
 [0, 0, 1, 0, 0, -1, 0]&&&
\\
(-1,17,-59),&
 [0, 0, 1, 0, 0, -1, 0]&&&
\\
(-1,19,-67),&
 [0, 0, 1, 0, 0, -1, 0]&&&
\\
(-1,22,-79),&
 [0, 0, 1, 0, 0, -1, 0]&&&
\\
(1,3,-3),&
 [0, -1, 1, 0, 0, -1, 0],&
 [0, 0, 0, 1, -1, 0, 1],&
 [0, 0, 1, 0, 0, -1, 1],&
[0, 1, 0, 0, 0, 0, 1], \\&&
 [0, 0, 0, 1, -1, 0, -1], &
 [1, -1, 2, -1, 1, -1, 0],&
 [1, 1, -2, -1, 1, 1, 0],\\&&
 [0, 1, -1, 0, 0, 0, 0],  &
[0, 0, 1, 0, 0, 0, 0]&
\\
(1,7,-19),&
 [0, -1, 1, 0, 0, -1, 0]&&&
\\
(1,11,-35),&
 [0, -1, 1, 0, 0, -1, 0]&&&
\\
(1,15,-51),&
 [0, -1, 1, 0, 0, -1, 0]&&&
\\
(1,19,-67),&
 [0, -1, 1, 0, 0, -1, 0]&&&
\end{array}
\]

\[
\begin{array}{lllll}
(3,5,-3),&
[0, -2, 1, 0, 0, -1, 0], &
[0, 1, -2, 0, 0, 1, -1], &
 [0, -1, -1, 0, 0, 0, -1], &
 [-1, 1, 0, 0, -1, 1, 0],  \\&&
 [1, 1, 0, 0, 1, 1, 0],&
 [0, 0, -2, 0, -1, 1, -1],&
 [0, 0, -2, 0, 1, 1, -1], \\&&
 [0, 1, -1, 0, 0, 0, 0],  &
 [0, 0, 1, 0, 0, 0, 0]&
\\
(3,6,-7),&
[0, -2, 1, 0, 0, -1, 0]&&&
\\
(3,7,-11),&
 [0, -2, 1, 0, 0, -1, 0]&&&
\\
(3,9,-19),&
 [0, -2, 1, 0, 0, -1, 0]&&&
\\
(3,14,-39),&
 [0, -2, 1, 0, 0, -1, 0]&&&
\\
(3,15,-43),&
 [0, -2, 1, 0, 0, -1, 0]&&&
\\
(3,18,-55),&
 [0, -2, 1, 0, 0, -1, 0]&&&
\\
(3,21,-67),&
 [0, -2, 1, 0, 0, -1, 0]&&&
\\
(3,25,-83),&
 [0, -2, 1, 0, 0, -1, 0]&&&
\\
(5,11,-11),&
 [0, -3, 1, 0, 0, -1, 0]&&&
\\
(5,19,-43),&
 [0, -3, 1, 0, 0, -1, 0]&&&
\\
(5,23,-59),&
 [0, -3, 1, 0, 0, -1, 0]&&&
\\
(7,15,-3),&
[0, -4, 1, 0, 0, -1, 0],&&&
\\
(7,17,-11),&
[0, -4, 1, 0, 0, -1, 0],&&&
\\
(7,18,-15),&
[0, -4, 1, 0, 0, -1, 0],&&&
\\
(7,19,-19),&
[0, -4, 1, 0, 0, -1, 0],&&&
\\
(7,22,-31),&
[0, -4, 1, 0, 0, -1, 0],&&&
\\
(7,23,-35),&
[0, -4, 1, 0, 0, -1, 0],&&&
\\
(7,25,-43),&
[0, -4, 1, 0, 0, -1, 0],&&&
\\
(9,23,-3),&
[0, -5, 1, 0, 0, -1, 0],&&&
\end{array}
\]

\section{Another solution}

In addition to  [0,1,0,0,0,0,0] 
in all cases considered, there appeared another solution with 
[0,$k$,1,0,0,-1,0] where $k$ seems to be related to $a$.
The above table was constructed to indicate this relation clearly.
The vector [0,$k$,1,0,0,-1,0] yields the element
\[
\alpha(k+\omega)-\alpha^3.
\]
If we try find the corresponding $X,Y,Z$ in the form (\ref{kpqpq}),
we get  $X=k+\omega=P^2-\delta Q^2,Y=0=PQ,Z=-1=Q^2$, which is
not possible. But if we consider the negative of this element, that is 
\[
\alpha(-k-\omega)-\alpha^3,
\]
then we have $X=-k-\omega=P^2-\delta Q^2,Y=0=PQ,Z=1=Q^2$
which has the solution $P=0,Q=1$, whence 
\[
-k-\omega=-\delta.
\]
This implies 
\[
k+\frac{1+\sqrt{m}}{2}=\frac{-a+\sqrt{m}}{2},
\]
implying 
\[
k=-\frac{a+1}{2}.
\]
Indeed, this is shown by the examples.
There remained to prove it formally.

\begin{theorem}
If $m=a^2-4b+8$ is square free, then
\begin{equation}
\gamma=\left(\frac{a+1}{2}-\omega\right)\alpha-\alpha^3
\label{gggh}
\end{equation}
generates a power integral basis in $K$.
\label{nnn}
\end{theorem}

\noindent
{\bf Proof.}
We have
\[
\frac{a+1}{2}-\omega=\frac{a-\sqrt{m}}{2}.
\]

Using the notation of Lemma \ref{relindex} we have
\[
\gamma^{(1,j)}=\frac{a-\sqrt{m}}{2}\cdot  \alpha^{(j)}+(\alpha^{(j)})^3,
\]
\[
\gamma^{(2,j)}=\frac{a+\sqrt{m}}{2}\cdot \overline{\alpha^{(j)}}
+(\overline{\alpha^{(j)}})^3,
\]
for $1\le j\le 4$.  As we have seen above, the element (\ref{gggh})
satisfies $I_{K/M}(\gamma)=1$ (it comes from a valid representation of 
$X,Y,Z$ by suitable $P,Q$ in (\ref{kpqpq})), therefore we only have to check 
$J(\gamma)=1$.
Using symmetric polynomials ($\alpha^{(j)}$ satisfies $x^4-\delta x^2+1=0$
and $\overline{\alpha^{(j)}}$ satsifies $x^4-\overline{\delta} x^2+1=0$)
we calculated 
\[
\prod_{j_1=1}^4\prod_{j_2=1}^4 |\gamma^{(1,j_1)}-\gamma^{(2,j_2)}|,
\]
and, making all possible simplifications, we found that it is 
equal to $m^2$. This calculation was also performed by Maple, and 
after several optimizations it took a negligible time. Note that 
this calculation of the above degree 16 polynomial must be
made very carefully, otherwise it results unusable formulas.
$\square$

\end{document}